\documentstyle[11pt]{article}

\newtheorem{theorem}{Theorem}[section]
\newtheorem{lemma}[theorem]{Lemma}
\newtheorem{corollary}[theorem]{Corollary}
\newtheorem{proposition}[theorem]{Proposition}

\newtheorem{rem}[theorem]{Remark}
\newenvironment{remark}{\begin{rem}\rm}{\end{rem}}

\newtheorem{define}[theorem]{Definition}
\newenvironment{definition}{\begin{define}\rm}{\end{define}}

\newtheorem{ex}[theorem]{Example}
\newenvironment{example}{\begin{ex}\rm}{\end{ex}}

\newenvironment{proof}{\par {\em Proof.}}{\hfill }

\newcommand{\Tr}{\mbox{Tr}} 
\newcommand{\tr}{\mbox{tr}}
\newcommand{\id}{\mbox{id}}

\newcommand{\Index}{\mbox{Index}\,}
\newcommand{\End}{\mbox{End}\,}

\newcommand{\eps}{\varepsilon}
\newcommand{\C}{\Bbb{C}}

\newcommand{\lact}{\triangleright}
\newcommand{\ract}{\triangleleft}
\newcommand{\la}{\langle}
\newcommand{\ra}{\rangle}
\newcommand{\flip}{\varsigma}
\newcommand{\nv}{\vec{m}}
\newcommand{\rtimes}{{>\!\!\!\triangleleft}}
\newcommand{\ltimes}{{\triangleright\!\!\!<}}
\newcommand{\smrtimes}{>\!\!\triangleleft}
\newcommand{\equivalent}{\Leftrightarrow}
\newcommand{\implies}{\Rightarrow}
\renewcommand{\L}{\Lambda}

\def\underset#1#2{\mathop{#2}\limits_{#1}}

\newcommand{\1}{_{(1)}}
\newcommand{\2}{_{(2)}}
\newcommand{\3}{_{(3)}}

\font\twlmsbm=msbm10 scaled \magstep1
\font\egtmsbm=msbm8
\font\sixmsbm=msbm6
\newfam\msbmfam
 
\textfont\msbmfam=\twlmsbm
\scriptfont\msbmfam=\egtmsbm   
\scriptscriptfont\msbmfam=\sixmsbm

\def\Bbb#1{{\fam\msbmfam\relax#1}}

\title{\bf Duality for actions of  weak Kac algebras \\
and  crossed product inclusions of  II${}_1$ factors }
\author{Dmitri Nikshych \thanks{UCLA, Department of Mathematics,
405 Hilgard Avenue, Los Angeles, CA 90095-1555;  
E-mail: nikshych@math.ucla.edu} }
\date{October 5, 1998}

\begin{document}
\maketitle


\begin{abstract}
We show that indecomposable weak Kac algebras are free over their
Cartan subalgebras and prove a duality theorem for their actions. 
Using this result, for any biconnected weak Kac algebra 
we construct a minimal action on the hyperfinite II${}_1$ factor. 
The corresponding crossed product
inclusion of II${}_1$ factors has depth 2 and an integer index.
Its first relative commutant is, in general, non-trivial, so
we derive some arithmetic properties of weak Kac algebras
from considering reduced subfactors.
 
\end{abstract}


\begin{section}
{Introduction}

It is well understood now that Kac algebras (Hopf $C^*$-algebras)
are closely related with the subfactors theory : it was announced by Ocneanu 
and proved in \cite{S},\cite{L},\cite{EN},\cite{D} that irreducible depth 2 
inclusions of type II factors come from crossed products with Kac algebras. 
This result was recently extended to the case of general (i.e., not necessarily
irreducible) finite index depth 2  subfactors in \cite{NV1}, where it was shown 
that if $N\subset M \subset M_1 \subset M_2\subset \dots$ is the Jones tower 
constructed from such a subfactor $N\subset M$, then $K=M^\prime \cap M_2$ 
has a natural structure of a finite-dimensional weak Kac algebra (if the index
$[M:N]$ is integer) or a weak Hopf $C^*$-algebra (if $[M:N]$ is non-integer) 
and there is a minimal action of $K$ on $M_1$ such that $M$ is the fixed point 
subalgebra of $M_1$ and $M_2$ is isomorphic to the crossed product of $M_1$ 
and $K$. This result establishes an injective correspondence 
between finite index depth 2 subfactors of a given II${}_1$ factor and weak 
$C^*$-Hopf algebras.

It is natural to ask if this correspondence is one-to-one
in the case of the hyperfinite II${}_1$ factor.
Note that in \cite{Y2}  Yamanouchi constructed
an outer action of any finite dimensional Kac algebra $K$
on the  hyperfinite II${}_1$ factor $R$, here the outerness means that 
$R'\cap R \rtimes K =\C$, i.e.\ the first relative commutant 
of the crossed product inclusion $R \subset R\rtimes K$ is minimal. 
His construction used the Takesaki duality for actions of
Kac algebras \cite{ES-Takesaki}.

In this work we extend this result to weak Kac algebras, i.e.,
we show that any weak Kac algebra has a minimal action on $R$.
Finite dimensional weak Kac algebras generalize both finite groupoid
algebras and  usual Kac algebras. Note that weak Kac algebra is a special 
case of a weak Hopf $C^*$-algebra introduced in \cite{BNSz} and \cite{Nill},
which is characterized by the property $S^2 =\id$.
It was shown in \cite{NV} that the category of weak Kac algebras
is equivalent to the categories of generalized Kac algebras of
T.~Yamanouchi \cite{Y1}  and Kac bimodules (an algebraic version of 
Hopf bimodules \cite{EV}). Compared with these objects,
the advantage of the language of weak Kac algebras is that 
their definition is transparently self-dual, so it is easy to work
with both weak Kac algebra and its dual simultaneously.
\medskip

The paper is organized as follows.

In Section 2 we collect the necessary definitions and
facts about weak Kac algebras, their actions and crossed-products,
and their Cartan subalgebras;
we also give a brief description of the basic construction for
$*$-algebras.

In Section 3 we introduce a $\lambda$-Markov condition for weak Kac
algebras. A weak Kac algebra $K$ satisfies the $\lambda$-Markov condition 
if the normalized Haar trace on $K$ is the $\lambda$-Markov trace
for the inclusion $K_s \subset K$, where $K_s$ is the source
Cartan subalgebra of $K$.
This condition is automatically satisfied if $K$ is indecomposable, i.e.,
not isomorphic to the direct sum of two weak Kac algebras.
Theorem~\ref{Markov criterion} shows that being $\lambda$-Markov
is equivalent to the freeness of $K$ over its Cartan subalgebras;
in particular, $\lambda^{-1}$ must be a positive integer.
As a corollary, we obtain that indecomposable weak Kac algebras of prime
dimension are group algebras of cyclic groups, which extends the
well-known result of Kac \cite{Kac}.

Also in this section we introduce and study basic properties of connected
and biconnected weak Kac algebras, i.e., those for which the inclusion
$K_s \subset K$ is connected (resp.\ inclusions $K_s \subset K$
and  $K_s^* \subset K^*$ are connected). The latest class of 
indecomposable weak Kac algebras is the most
important for the applications to subfactors in Section 5, 
so we describe a way of constructing biconnected weak Kac algebras
from two-sided crossed-products introduced in \cite{Nill}.

The central result of Section 4 is a duality theorem for actions
of weak Kac algebras. This theorem is an analogue of the well known duality 
results for actions of locally compact groups \cite{NT}, 
Kac algebras \cite{ES-Takesaki}, and Hopf algebras \cite{BM}. 
It states that if $K$ satisfies the $\lambda$-Markov condition and acts on a 
$C^*$-algebra (von Neumann algebra) $A$, then the dual crossed product algebra
$(A\rtimes K)\rtimes K^*$ is isomorphic to $A \otimes M_{\lambda^{-1}}(\C)$.
Let us note that a similar result for depth 2 inclusions of von Neumann
algebras was proved in \cite{EV}.

In Section 5 for any  biconnected weak Kac algebra $K$ we construct a 
minimal action on the  hyperfinite II${}_1$ factor $R$ 
(where the minimality means that the relative commutant $R'\cap R\rtimes K$ 
is minimal). The resulting crossed product inclusion $R\subset R\rtimes K$ 
of II${}_1$ factors has depth 2 and an integer index $\lambda^{-1}$. 
We compute the standard invariant of this inclusion, an show,
in particular, that the first relative commutant is isomorphic to the
Cartan subalgebra of $K$ : $R'\cap R\rtimes K\cong K_s$.

Finally, in Section 6 we construct examples of irreducible subfactors reducing
the inclusion $R\subset R\rtimes K$ by the minimal projection in 
$K_s = R'\cap R\rtimes K$. In this way we can associate
an irreducible finite depth subfactor of $R$ with every irreducible 
representation of $K$ or $K_s$. This allows us to derive certain arithmetic 
properties of biconnected weak Kac algebras.

{\bf Acknowledgement.} The author is deeply grateful to L.~Vainerman
for many important discussions. He also would like to thank S.~Popa
and E.~Effros for helpful advices, and E.~Vaysleb for his useful
comments on the preliminary version of this paper.

\end{section}


\begin{section}
{Preliminaries}

\medskip
{\bf 2.1 Weak Kac algebras \cite{BNSz}, \cite{NV}.}
\medskip

Throughout this paper all weak Kac algebras are supposed to be
finite-dimensional.

The notion of weak Kac algebra \cite{NV} is a special case of a more general 
concept of weak $C^*$-Hopf algebra introduced in \cite{BNSz}; see \cite{NV}
for a discussion on equivalence of weak Kac algebras with other
algebraic versions of quantum groupoid ( generalized Kac algebras of 
T.~Yamanouchi \cite{Y1} and Kac bimodules).
\medskip

A {\em weak Kac algebra} $K$ is a finite dimensional $C^*$-algebra equipped
with the following linear maps
\begin{enumerate}
\item[(a)] {\em comultiplication} $\Delta : K\to K \otimes K$,
\item[(b)] {\em counit} $\eps : K \to \C$,
\item[(c)] {\em antipode} $S : K\to K$,
\end{enumerate}
where $\Delta$ is a (not necessarily unital) homomorphism of $C^*$-algebras,
$\eps$ is a positive (not necessarily multiplicative) functional, 
$S$ is a $*$-preserving anti-multiplicative and anti-comultiplicative 
involution (i.e., $S^2=\id$) such that the following identities hold 
(we denote $\eps_s(x) = (\id\otimes \eps)((1\otimes x)\Delta(1) )$ and
$\eps_t(x) = (\eps\otimes \id)(\Delta(1) (x\otimes 1) )$ ) :
\begin{enumerate}
\item[(1)] $(\Delta\otimes \id)\Delta = (\id\otimes \Delta)\Delta, \quad
            (\eps\otimes\id)\Delta =\id = (\id\otimes\eps)\Delta$,
           \newline i.e., $K$ is a coalgebra,
\item[(2)] $\eps_s(x)y = (\id\otimes \eps)((1\otimes x)\Delta(y) )$,
\item[(3)] $(\eps_s\otimes\id) \Delta(x) = (1\otimes x)\Delta(1)$,
\item[(4)] $m(S\otimes\id)\Delta(x) = \eps_s(x),  \quad x,\,y\in K,$
\end{enumerate}
where $m$ denotes multiplication. 

The following identities are equivalent to the above axioms (2)--(4) 
respectively.
\begin{enumerate}
\item[(2')] $x\eps_t(y) =  (\eps\otimes \id)(\Delta(x) (y\otimes 1) )$,
\item[(3')] $(\id\otimes\eps_t) \Delta(x) = \Delta(1)(x\otimes 1)$,
\item[(4')] $m(\id\otimes S)\Delta(x) = \eps_t(x), \quad x,\,y\in K.$
\end{enumerate}
\medskip

The dual vector space $K^*$ has a natural structure of a weak Kac algebra
given by dualizing the structure operations of $K$ :
\begin{eqnarray*}
& & \la \phi\psi,\,x\ra = \la \phi\otimes\psi,\,\Delta(x) \ra 
\quad (multiplication), \\
& & \la \Delta(\phi),\,x \otimes y \ra =  \la \phi,\,xy \ra
\quad (comultiplication), \\
& & \la S(\phi),\,x\ra = \la \phi,\,S(x) \ra
\quad (antipode), \\
& &  \la \phi^*,\,x\ra = \overline{ \la \phi,\,S(x^*) \ra }
\quad (*-operation),
\end{eqnarray*}
for all $\phi,\psi \in K^*,\, x,y\in K.$ Unit is given by $\eps$
and counit is by $\phi \mapsto \la\phi,\, 1\ra$. 
\medskip

Below we collect the most important results of the theory of weak 
Kac algebras.  The proofs can be found in \cite{NV}.

The maps $\eps_s$ and $\eps_t$ are called {\em source} and {\em target}
{\em counital maps} respectively, we have  $\eps_s^2 = \eps_s$ and  
$\eps_t^2 = \eps_t$. Their images are unital $C^*$-subalgebras,
called {\em Cartan subalgebras} of $K$:
\begin{eqnarray*}
K_s &=& \{x\in K \mid \eps_s(x) =x \} 
= \{x\in K \mid \Delta(x) = 1\1 \otimes x 1\2 = 1\1 \otimes 1\2 x \},\\
K_t &=& \{x\in K \mid \eps_t(x) =x \} 
  =  \{x\in K \mid \Delta(x) = x 1\1 \otimes 1\2 = 1\1 x\otimes 1\2 \}.
\end{eqnarray*}

Cartan subalgebras commute : $[K_s,\, K_t] =0$; also we have 
$S\circ\eps_s = \eps_t\circ S$ and $S(K_s) =K_t$.
\medskip

Like usual finite-dimensional Kac algebras ( $=$ Hopf $C^*$-algebras),
weak Kac algebras have integrals in the following sense.

\begin{enumerate}
\item[]
There exists a unique projection $p_\eps\in K$,
called a {\em Haar projection}, such that for all $x\in K$ :
$$
p_\eps x = p_\eps \eps_s(x),\quad 
x p_\eps  = \eps_t(x)p_\eps ,\quad 
\eps_s(p_\eps )=\eps_t(p_\eps )=1. 
$$
\item[]
There exists a unique fathfull trace $\tau$ on $K$, 
called a {\em normalized Haar trace}, such that
$$
(\tau\otimes \id)\Delta = (\tau\otimes \eps_s)\Delta,\quad
(\id\otimes \tau)\Delta = (\eps_t\otimes \tau)\Delta,\quad
\tau\circ \eps_s = \tau\circ \eps_t = \eps.
$$
\end{enumerate}

Note that normalized Haar projection and trace are unimodular, i.e. 
$S(p_\eps ) = p_\eps$ and $\tau\circ S =\tau$.
By duality,  $\tau$ is the normalized Haar projection
for the dual weak Kac algebra $K^*$.

The maps
\begin{eqnarray*}
E_s &:& K \to K_s : E_s(x) = (\tau\otimes \id)\Delta(x), \\
E_t &:& K \to K_t : E_t(x) = (\id\otimes \tau)\Delta(x) 
\end{eqnarray*}
define $\tau$-preserving conditional expectations 
(see 2.3.\ for the definition) from $K$ to Cartan subalgebras.
\medskip

To fix the notation in what follows, let
$$
K\cong \oplus_{i=1}^N\,M_{d_i}(\C), \qquad
K_s \cong K_t \cong \oplus_{\alpha=1}^L\,M_{m_\alpha}(\C),
$$
and let $\{e_{kl}^{(i)}\}$ $(i=1\dots N;\,  k,l=1,\dots, d_i)$ be a system of
matrix units in $K$, $\{f_{rs}^{(\alpha)}\}$ in $K_s$,
and  $\{g_{rs}^{(\alpha)}\}$ in $K_t$
$(\alpha =1,\dots, L;\,  r,s = 1,\dots, m_\alpha)$. By \cite{NV} we have~:
\begin{eqnarray*}
\Delta(p_\eps)
&=& \sum_i\,\frac{1}{d_i}\,\sum_{kl}\,e_{kl}^{(i)} \otimes S(e_{lk}^{(i)}), \\
\Delta(1)
&=& \sum_\alpha\,\frac{1}{m_\alpha}\,\sum_{rs}\,f_{rs}^{(\alpha)}
\otimes S(f_{sr}^{(\alpha)})
=  \sum_\alpha\,\frac{1}{m_\alpha}\,\sum_{rs}\,S(g_{sr}^{(\alpha)})
\otimes g_{rs}^{(\alpha)}.
\end{eqnarray*}
In particular, $p_\eps$ is  cocommutative, i.e., 
$\Delta(p_\eps) = \flip \Delta(p_\eps)$, where $\flip$ is the flip
on the tensor product $K\otimes K$. 

Also we denote $\Lambda = (\Lambda_{ij})$ the $(L\times N)$ inclusion matrix 
\cite{GHJ} of $K_s$ (or $K_t$) into $K$.
\bigskip

{\bf  2.2 Actions, dual actions,  and crossed-products \cite{NSzW}.}
\medskip

By a $*$-algebra we understand an associative algebra over $\C$
equipped with an antilinear anti-automorphism of order $2$ (involution),
$x\mapsto x^*$.

Action of a weak $C^*$-Hopf algebra on a $*$-algebra
was defined in \cite{NSzW}. We slightly modify that definition, since
we consider only those actions for which the map $x \mapsto (x\lact 1)$
($x \mapsto (1 \ract x)$)  is injective on a Cartan subalgebra. 
We need definitions of left and right actions.
\medskip

\label{action}
A {\em left} (resp. {\em right}) {\em action} of a weak Kac algebra $K$ on a 
unital $*$-algebra $A$ is a linear map 
$$
K\otimes A \ni h\otimes a \mapsto (h\lact a)\in A 
\quad (\mbox{resp. } A\otimes K \ni  a\otimes h \mapsto (a \ract h) \in A),
$$
defining a structure of a left (resp. right) $K$-module on $A$ such that :
\begin{enumerate}
\item[(1)]
$h\lact ab = (h\1\lact a)(h\2\lact b) \quad$
 (resp. $ab \ract h = (a \ract h\1)(b \ract h\2)$),
\item[(2)]
$(h\lact a)^* = Sh^*\lact a^* \quad $
 (resp. $(a \ract h)^* = a^* \ract Sh^*$),
\item[(3)]
$h\lact 1 = \eps_t(h)\lact 1$, and  $h\lact 1 =0$ iff $\eps_t(h)=0$ \newline
 (resp. $1 \ract h = 1\ract \eps_s(h)$, and $1 \ract h = 0$ iff $\eps_s(h)=0$). 
\end{enumerate}
\medskip
 
If $A$ is a $C^*$-algebra or a von Neumann algebra then
we also require that the map $a \mapsto  (h\lact a)$ (resp.
$a \mapsto (a \ract h)$) to be norm continuous or weakly continuous 
for all $h\in K$.

Note that the map $z\mapsto (z\lact 1)$  (resp. $z \mapsto (1 \ract z)$)
defines an injective  $*$-homomorphism from $K_t$ (resp. $K_s$) to $A$. 
Thus, $A$  must contain a $*$-subalgebra isomorphic to a Cartan 
subalgebra of $K$.
\newline
A {\em trivial}  left (resp. right) action of $K$ on $K_t$ (resp. $K_s$) 
is given by 
$$
h\lact a = \eps_t(ha) \quad ( \mbox{resp. } a \ract h = \eps_s(ah)), 
\qquad h\in K, \,a\in K_t\,(\mbox{resp. } K_s). 
$$ 
A {\em dual} left (resp. right) action of $K^*$ on $K$  is given by
$$
\phi\lact h = h\1 \la \phi,\,h\2 \ra  
\quad ( \mbox{resp. }  h \ract \phi = \la \phi,\,h\1 \ra h\2),
\qquad \phi\in K^*,\,h\in K.
$$

Given a left (resp. right) action of $K$ on a $*$-algebra $A$,
there is a left (resp. right) {\em crossed product} $*$-algebra
$A\rtimes K$ (resp. $K\ltimes A$) constructed as follows. 
As a $\C$-vector space it is $A \otimes_{K_t} K$ (resp. $K \otimes_{K_s} A$),
where  $K$ is a  left (resp. right) $K_t$-module  (resp. $K_s$-module) 
via multiplication and $A$ is a right $K_t$-module (resp. left
$K_s$-module) via multiplication by image of $K_t$ (resp. $K_s$)  under
$z\mapsto (z\lact 1)$  (resp. $z \mapsto (1 \ract z)$); that is,
we identify 
$$
a(z\lact 1) \otimes h  \equiv  a \otimes zh  \quad 
( \mbox{resp. } hz \otimes a  \equiv h \otimes (1\ract z)a ),
$$
for all $a\in A,\,h\in K,\,z\in K_t\,(K_s)$.
Let $[a\otimes h]$ (resp. $[h\otimes a]$) denote the class of $a\otimes h$
(resp. $h\otimes a$). A $*$-algebra structure is defined by
\begin{eqnarray*}
[a\otimes h][b\otimes k] &=& [a(h\1\lact b) \otimes h\2 k], \quad
[a\otimes h]^* = [(h\1^*\lact a^*) \otimes h\2^*] \\
( \mbox{resp. }
[h\otimes a][k\otimes b] &=& [hk\1 \otimes (a\ract k\2)b], \quad
[h\otimes a]^* = [h\1^* \otimes (a^* \ract  h\2^*)] ),
\end{eqnarray*}
for all $a,b\in A,\,h,k \in K$.
The maps $i_A: a \mapsto [a\otimes 1_K]$ (resp. $a \mapsto [1_K\otimes a]$)
and  $i_K : h \mapsto [1_A\otimes h]$  (resp. $h \mapsto [ h\otimes 1_A]$)
are inclusions of $*$-algebras such that $A\rtimes K = i_A(A)\,i_K(K)$
(resp. $K \ltimes A = i_K(K)\,i_A(A)$).
Moreover, if $A$ is a $C^*$-algebra (von Neumann algebra), 
then the crossed product  is naturally $*$-isomorphic to a norm closed (weakly
closed) $*$-algebra of operators on some Hilbert space, i.e., it becomes
a  $C^*$-algebra (von Neumann algebra).

For the crossed products constructed from the trivial
actions of $K$ on Cartan subalgebras we have 
$$
K_t\rtimes K \cong K \quad \mbox{and} \quad K \ltimes K_s \cong K.
$$ 
A left (resp. right) {\em dual action} of $K^*$ on the crossed product 
$A\rtimes K$  (resp. $K \ltimes A$) is defined as
$$
\phi\lact [a\otimes h] = [a \otimes (\phi\lact h)] \quad
( \mbox{resp. } [h\otimes a] \ract \phi = [(h\ract \phi) \otimes a] ),
$$ 
for all $a\in A,\,h\in K,\, \phi\in K^*$.
The action of $K^*$ on $K$ defined above is dual to the trivial action 
of $K$ on a Cartan subalgebra.
\bigskip


{\bf 2.3 The basic construction for $*$-algebras \cite{W}.}
\medskip

Let $B$ be a unital $*$-algebra, $A$ be its $*$-subalgebra containing
the unit of $B$. A {\em conditional expectation} $E : B \to A$ is a faithful
(i.e, such that $E(Bb)=0$ implies $b=0$, for $b\in B$ )
linear *-preserving map satisfying
$$
E(ab) = a E(b),\quad E(ba)=E(b)a,\quad \mbox{and} \quad E(a) =a,
$$
for all $a\in A,\, b\in B$. 
A finite family $\{\,u_1,\dots, u_n\,\}\subset B$ is called a 
{\em quasi-basis} for $E$ if
$$
b  = \Sigma_i\,u_i E(u_i^* b) \qquad \mbox{for all } b\in B.
$$
It is called a {\em basis} if ``coefficients'' $E({u_i^* b})$
are unique, i.e. if $\sum_i\, u_ia_i =0,\,a_i\in A
\Leftrightarrow a_i =0\,(\forall i)$. A conditional expectation $E$ is of
{\em finite-index type} if there exists a quasi-basis for $E$. In this case
{\em index} of $E$ is defined as
$$
\Index E = \Sigma_i\,u_i u_i^* \in B.
$$ 
$\Index E$ belongs to the center of $B$
and does not depend on the choice of quasi-basis. 

The {\em basic construction} for $E$ is a $*$-algebra $B\otimes_A B$
with multiplication and involution  given by
$$
(b_1\otimes b_2)(b_3\otimes b_4) = b_1E(b_2b_3) \otimes b_4, \quad
(b_1\otimes b_2)^* = b_2^*\otimes b_1^*,
$$
for all $b_1,b_2,b_3,b_4$ in $B$. Note that the unit of this algebra
is $\Sigma_i u_i\otimes u_i^*$, where $\{ u_i\}$ is the quasi basis
for $E$.
\medskip 

In what follows we consider only conditional expectations of 
finite-index type for which a basis (not just a quasi-basis) exists.

In this case  $B\otimes_A B$ is canonically isomorphic to $\End_A^r(B)$,
the algebra of endomorphisms of $B$ viewed as a right $A$-module;
the isomorphism $\phi : B\otimes_A B \to \End_A^r(B)$ is given by
$$
\phi(b_1\otimes b_2)(b) = b_1 E(b_2 b), \qquad b,b_1,b_2\in B.
$$
$B$ is canonically identified with the subalgebra of left multiplication 
operators ($x\mapsto bx$ for $x\in B$) in $\End_A^r(B)$.
Clearly, $\End_A^r(B) \cong M_n(A)$, the $*$-algebra of $(n\times n)$
matricies over $A$, since $B$ is  free of rank $n$ over $A$.

Note that $e_A = E \in \End_A^r(B)$ is a projection such that 
\begin{enumerate}
\item[(1)]
$e_A b e_A = E(b) e_A$ for all $b\in B$,
\item[(2)]
the map $A\ni a \mapsto ae_A\in \End_A^r(B)$ is injective,
\end{enumerate}
and $\End_A^r(B)$ is generated by $B$ and $e_A$.  

Conversely, if $C$ is a $*$-algebra containing $B$ as a unital
$*$-subalgebra and generated by $B$ and some projection $e_A$
satisfying properties (1) and (2) above then $C$ is canonically
$*$-isomorphic to $B\otimes_A B$, the isomorphism is given by
$b_1e_Ab_2 \mapsto b_1\otimes b_2$.

Due to this fact, we will denote the basic construction for $E$
by  $\la B,\,e_A\ra$.

When $A$ and $B$ are $C^*$ -algebras (von Neumann algebras) then
$\la B,\,e_A\ra$ naturally becomes a $C^*$-algebra (von Neumann algebra).
\end{section}


\begin{section}
{$\lambda$-Markov Condition and Connected Weak Kac Algebras}

We show that indecomposable weak Kac algebras have the Markov property
(i.e., $E_s(p_\eps)$ is a scalar). We prove that this property is equivalent
to existence of a basis for $E_s$, i.e., to freeness of $K$ over the Cartan
subalgebra $K_s$, and can be expressed  in terms of the inclusion matrix
of $K_s \subset K$.

We also introduce notions of connected and biconnected weak Kac algebras
which are important for applications to the theory of subfactors.

\begin{definition}
\label{indecomposability}
A weak Kac algebra $K$ is {\em decomposable} if it is isomorphic to the direct
sum of two weak Kac algebras, $K\cong K_1 \oplus K_2$; otherwise $K$ is
{\em indecomposable}.
\end{definition}

Clearly, decomposable weak Kac algebras are of little interest
in themselves. Indecomposability can be expressed in terms of the algebra
$K_s \cap K_t \cap Z(K)$, the {\em hypercenter} of $K$ \cite{NSzW} (here
$Z(K)$ is the center of $K$).

\begin{proposition}
\label{hypercenter}
$K$ is indecomposable iff the hypercenter of $K$ is trivial, i.e.
$K_s \cap K_t \cap Z(K) = \C$.
\end{proposition}
\begin{proof}
If $q$ is a projection in $K_s \cap K_t \cap Z(K)$, then $S(q) =q$,
$\Delta(q) = (q\otimes q)\Delta(1)$, and $\Delta(qh) = (q\otimes q)\Delta(h)$
for all $h\in K$. Therefore, $qK$ and $(1-q)K$ are weak Kac algebras and
$K = qK \oplus (1-q)K$. Conversely, if $K$ is decomposable and $K_1$
is its direct summand, then $\eps_s(1_{K_1}) = \eps_t(1_{K_1}) = 1_{K_1}$
and the unit of $K_1$ belongs to the hypercenter of $K$.
\end{proof}

It turns out that indecomposable weak Kac algebras satisfy a crucial
$\lambda$-Markov condition.

\begin{definition}
A weak Kac algebra $K$ satisfies a $\lambda$-Markov
condition if  
$$
E_s(p_\eps) =  E_t(p_\eps) = \lambda
$$
for some $\lambda >0$ (note that since $p_\eps$ is cocommutative,
we always have $E_s(p_\eps) =  E_t(p_\eps)$ ). 
\end{definition}

\begin{proposition}
\label{Markovianity}
If $K$ is indecomposable, then it satisfies the $\lambda$-Markov    
condition for some $\lambda$.
\end{proposition}
\begin{proof}
Note that 
$$
E_s(p_\eps) =  E_t(p_\eps) =
\sum_i\,\frac{1}{d_i}\sum_{k}\,e_{kk}^{(i)} \tau(e_{kk}^{(i)}) =
\sum_i\,\frac{\tau_i}{d_i}\sum_{k}\,e_{kk}^{(i)} \in Z(K),
$$
where $\tau_i =  \tau(e_{kk}^{(i)})$ does not depend on $k$.
Therefore, if $E_s(p_\eps) \neq \lambda$, then the hypercenter is
non-trivial and $K$ is decomposable by Proposition~\ref{hypercenter}.
\end{proof} 

The following theorem describes the  $\lambda$-Markov condition
in several equivalent ways.

\begin{theorem}
\label{Markov criterion}
The following conditions are equivalent :
\begin{enumerate}
\item[(i)]
K satisfies the $\lambda$-Markov condition,
\item[(ii)] 
$\tau = \lambda\,\Tr$ where $\Tr$ is the trace of the left 
regular representation of $K$ on itself,
\item[(iii)]
$(\L\L^t)\nv =\lambda\,\nv$, where $\nv = (m_1\dots m_L)$
is the dimension vector of a Cartan subalgebra, and $\L$ is
the $L\times M$ inclusion matrix of $K_s\subset K$.
\item[(iv)]
$n=\lambda^{-1}$ is an integer and there is a basis 
$\{x_\nu\}_{\nu=1\dots n}$
for $E_s$, i.e a basis of $K$ over $K_s$  such that 
$x =\sum_\nu\,x_\nu E_s(x_\nu^* x)$  for all $x\in K$,
 \item[(v)]
$n=\lambda^{-1}$ is an integer and there is a basis 
$\{y_\nu\}_{\nu=1\dots n}$
for $E_t$, i.e a basis of $K$ over $K_t$  such that 
$y =\sum_\nu\, y_\nu E_t(y_\nu^* y)$  for all $y\in K$.
\end{enumerate}
\end{theorem}
\begin{proof}
(i)$\equivalent$(ii). 
As we have seen in the proof of Proposition~\ref{Markovianity},
$E_s(p_\eps) = \lambda$ iff there is $\lambda$ such that 
$\tau(e_{kk}^{(i)}) = \lambda\,d_i$, i.e., $\tau = \lambda\,\Tr$.

(ii)$\equivalent$(iii).
It suffices to prove that (iii) holds true if and only if $\Tr$ is normalized 
by conditions $(\Tr\otimes\id)\Delta(1) = (\id\otimes\Tr)\Delta(1) 
=\lambda^{-1}$ (it was shown in \cite{NV} that 
$(\Tr\otimes \id)\Delta = (\Tr\otimes \eps_s)\Delta$ and  
$(\id\otimes \Tr)\Delta = (\eps_t\otimes \Tr)\Delta$).
Since $\Tr\circ S = \Tr$, we have
$$
(\Tr\otimes\id)\Delta(1) = \sum_{\alpha=1}^K\,\frac{1}{m_\alpha}\,
\sum_{r}\, \Tr(g_{rr}^{(\alpha)}) g_{rr}^{(\alpha)}
= \sum_{\alpha=1}^K\,\frac{1}{m_\alpha}
( \sum_i\,\L_{\alpha i}d_i ) \sum_{r}\,g_{rr}^{(\alpha)}
$$
$$
(\id\otimes\Tr)\Delta(1) =  \sum_{\alpha=1}^K\,\frac{1}{m_\alpha}\,
\sum_{r}\, f_{rr}^{(\alpha)} \Tr(f_{rr}^{(\alpha)}) 
= \sum_{\alpha=1}^K\,\frac{1}{m_\alpha}
( \sum_i\,\L_{\alpha i}d_i )  \sum_{r}\,f_{rr}^{(\alpha)}.
$$
This shows that (ii) is equivalent to the following condition :
$$
\Sigma_i\,\L_{\alpha i}d_i =\lambda^{-1}\, m_\alpha,\quad 
\alpha=1,\dots, L.
$$
But $d_i = \sum_{\beta =1}^L \,\L_{\beta i} m_\beta$ since
the inclusion $K_s \subset K$ is unital. Hence, we can rewrite
the last condition as 
$$
\sum_{i=1}^N\, \sum_{\beta =1}^L\, 
\L_{\alpha i}\, \L_{\beta i}\, m_\beta =\lambda^{-1}\, m_\alpha,
\quad \alpha = 1\dots L,
$$
which means precisely that $(\L\L^t)\nv =\lambda^{-1}\nv$.

(iii)$\implies$(iv).
It is clear that $\lambda^{-1}$ is a positive rational number since all entries
of $(\L\L^t)$ and $\nv$ are positive integers. On the other hand, 
$\lambda^{-1}$ is an algebraic integer, since it is an eigenvector 
of the integer matrix $(\L\L^t)$, therefore, $\lambda^{-1}$ is an integer.

For all $\alpha= 1,\dots, L$ and $r=1,\dots, m_\alpha$ define 
$K_{\alpha r} = K f_{rr}^{(\alpha)}$. We have 
$$
\dim(K_{\alpha r}) = \Tr(f_{rr}^{(\alpha)}) 
=\sum_i\,a_{\alpha i}d_i = n  m_\alpha.
$$
For all $y,\, z \in K_{\alpha r}$ :
$$
E_s(y^*z) = f_{rr}^{(\alpha)} E_s(y^*z) f_{rr}^{(\alpha)}
= (y,\,z)\,f_{rr}^{(\alpha)},
$$
where $(y,\,z)$ is a scalar since $f_{rr}^{(\alpha)}$ is minimal in
$K_s$. Clearly, $(\cdot,\cdot)$ defines an inner product
in $K_{\alpha r}$, which is non-degenerate since $E_s$ is faithful.
Let us choose an orthonormal bases 
$\{ x^{\alpha r}_\mu \},\, (\mu=1,\dots,  n  m_\alpha)$
in  $K_{\alpha r}$, $\alpha=1,\dots L,\,r=1,\dots, m_\alpha$ in such a way
that 
$$
x^{\alpha t}_\mu = x^{\alpha r}_\mu f_{rt}^{(\alpha)} \quad
\mbox{for all }\quad
t,r=1\dots m_\alpha,\,\, \mu=1\dots  n  m_\alpha.
$$
Then we have the following relation
$$
E_s( (x^{\alpha r}_\mu)^* x^{\alpha' r'}_{\mu'} ) =
\delta_{\alpha\alpha'}\,\delta_{\mu\mu'} f_{rr'}^{(\alpha)}
\quad \mbox{for all } \quad \alpha,\alpha';\,\, \mu,\mu';\,\, r,r'.
$$
We claim that
$$
x_\nu = \sum_\alpha\,\sum_{rs}\,\frac{1}{\sqrt{m_\alpha}} \,
\exp \left(\frac{2sr \pi}{m_\alpha} i\right)\, x_{\nu+(s-1)n}^{\alpha r}
\qquad \nu=1,\dots, n,
$$ 
is a basis of $K$ over $K_s$. Indeed :
\begin{eqnarray*}
\sum\nolimits_\nu\, x_\nu E_s(x_\nu^* x_\mu^{\beta t}) 
&=& \sum\nolimits_\nu\,\sum\nolimits_{\alpha r s}\,
\frac{1}{m_\alpha}  x_{\nu+(s-1) n }^{\alpha r}
E_s( (x_{\nu+(s-1) n }^{\alpha r})^*  x_\mu^{\beta t}) \\
&=& \frac{1}{m_\beta}\,\sum\nolimits_r\, x_\mu^{\beta r} f_{rt}^{(\beta)} 
= x_\mu^{\beta t} 
\end{eqnarray*}
for all  $\beta=1\dots K,\, t=1\dots m_\beta,\, \mu=1\dots  n  m_\beta$.
Next,
\begin{eqnarray*}
E_s(x_\nu^* x_\kappa) 
&=&  \sum_{\alpha rr'ss'}\,\frac{1}{m_\alpha}
\exp\left(\frac{2(sr - s'r')\pi}{m_\alpha}\,i \right) \,
E_s( (x_{\nu+(s-1)n }^{\alpha r})^* x_{\kappa+(s'-1) n }^{\alpha r'} ) \\
&=& \delta_{\nu\kappa}  \sum_{\alpha rr's}\,\frac{1}{m_\alpha}
\exp\left( \frac{2s(r -r')\pi}{m_\alpha}\,i \right)\,  f_{rr'}^{(\alpha)} \\
&=& \delta_{\nu\kappa} \sum\nolimits_{\alpha r}\, f_{rr}^{(\alpha)}
= \delta_{\nu\kappa}.
\end{eqnarray*}
Since `$E_s$ -orthogonality' implies linear independence over $K_s$,
we conclude that $\{x_\nu \}$ is a basis for $E_s$.

(iv)$\implies$(iii).
If there is a basis for $E_s: K\to K_s$ then the basic construction 
$\la K,\, e_{K_s} \ra$ is isomorphic to $M_n(K_s)$.
This means that that the inclusion matrix $B$ of the inclusion
$K_s \subset \la K,\, e_{K_s} \ra$ satisfies $B\nv =\lambda^{-1} \nv$.
But  $B=\L\L^t$ \cite{GHJ}.

(iv)$\equivalent$(v).
We will prove (iv)$\implies$(v), the converse implication is completely
analogous. If $x =\sum_\nu\, x_\nu E_s(x_\nu^* x)$ then
$Sx^* = \sum_\nu\, Sx_\nu^* E_t(Sx_\nu Sx^*)$, 
since $E_t = S\circ E_s \circ S$, and we can take $y_\nu =Sx_\nu^*$,
$\nu = 1\dots \lambda$ as a basis for $E_t$.
\end{proof}

\begin{corollary}
\label{tau is a Markov trace}
If the equivalent conditions of Theorem~\ref{Markov criterion}
are satisfied then $\tau$ is a $\lambda$-Markov trace for the inclusion
$K_t \subset K\,(K_s \subset K)$.
\end{corollary}
\begin{proof}
We need to show that $\L^t\L\vec{t} = \lambda^{-1}\,\vec{t}$, where
$\vec{t}$ is the `trace-vector' corresponding to $\tau$ (\cite{JS},
3.2.3(ii)). Since $\tau =\lambda\Tr$, we have $\vec{t} = \lambda\vec{d}$,
where $\vec{d} =(d_1\dots d_N)$ is the `dimension-vector' of $K$.
Using Theorem~\ref{Markov criterion}(iii) we compute
$$
\L^t\L\vec{t} = \lambda \L^t\L\vec{d} = \lambda  \L^t\L\L^t\vec{m}
= \L^t\vec{m} = \vec{d} =  \lambda^{-1}\,\vec{t}.
$$
\end{proof}

\begin{remark}
\label{properties of lambda}
\begin{enumerate}
\item[(i)]
Proposition~\ref{hypercenter} says that
$K$ is indecomposable iff the matrix $\L$ is indecomposable in the
sense of \cite{GHJ}. In this case
Theorem~\ref{Markov criterion}(iii) implies that $\nv$ is
the Perron-Frobenius eigenvector of the matrix  $(\L\L^t)$. It is well-known
that in this case the corresponding eigenvalue $\lambda^{-1}$ is equal to the 
spectral radius of $(\L\L^t)$, so 
$$
\lambda^{-1} =\Vert \L\L^t \Vert = \Vert \L \Vert^2. 
$$
\item[(ii)]
Theorem~\ref{Markov criterion}(iv) and (v) show that an indecomposable
weak Kac algebra $K$  is free over its Cartan subalgebras $K_s$ and $K_t$.
In particular, $\dim K_s$ divides $\dim K$ and
$$
\lambda^{-1} = \frac{\dim K}{\dim K_s}.
$$
\item[(iii)]
Conditional expectations 
$E_s$ and $E_t$ are of index-finite type and their index is an
integer scalar : $\Index E_s = \Index E_t = \lambda^{-1}$.
\end{enumerate}
\end{remark}

\begin{corollary}
\label{dim p}
If $K$ is indecomposable and $\dim K = p$, where $p$ is a prime,
then $K\cong \C\Bbb{Z}_p$, a group algebra of a simple abelian group.
\end{corollary}
\begin{proof}
Remark~\ref{properties of lambda}(ii) implies that Cartan subalgebras 
of $K$ must be 1-dimensional, so $K$ is a Kac algebra.  But in this case 
the result is well known \cite{Kac}.
\end{proof}
\medskip

The $\lambda$-Markov condition is invariant under duality.
\begin{proposition}
$K$ satisfies the $\lambda$-Markov condition iff $K^*$ satisfies
the $\lambda$-Markov  condition (with the same $\lambda$).
\end{proposition}
\begin{proof} 
Since $K$ satisfies the $\lambda$-Markov condition iff every its
indecomposable component does, it sufficed to prove this statement
in the case when $K$ is indecomposable. But this is trivial by
Proposition~\ref{Markovianity} and Remark~\ref{properties of lambda}(ii),
since $\dim K_s = \dim K_s^*$.
\end{proof}
\medskip

Connected weak Kac algebras (i.e., those with connected Bratteli
diagram of the inclusion $K_s\subset K$) form a subclass of indecomposable
weak Kac algebras important for the applications to
subfactors in section 5.

\begin{definition}
A  weak Kac algebra $K$ is {\em connected} if the inclusion
$K_s \subset K$ is connected, i.e., $K_s \cap Z(K) = \C$
(or, equivalently, $K_t\cap Z(K) = \C$), where $Z(\cdot)$ denotes
the center of an algebra.
$K$ is {\em biconnected} if both $K$ and $K^*$ are connected.
\end{definition}

\begin{proposition}[cf. \cite{NSzW}]
\label{connectedness}
The following conditions are equivalent :
\begin{enumerate}
\item[(i)]
$K$ is connected,
\item[(ii)]
$K_s^* \cap K_t^* = \C$,
\item[(iii)]
$p_\eps$ is a minimal projection in $K$ (i.e the counital representation
of $K$ (\cite{NV}, Section 2.2) is irreducible).
\end{enumerate}
\end{proposition}
\begin{proof}
(i)$\implies$(ii).
Suppose that there is $\beta\in K_s^* \cap K_t^*$,
$\beta\not\in \C$. Since the Cartan subalgebras commute, $\beta$
must belong to $Z(K_s^*)$, the center of $K_s^*$. Consider the element
$b\in K \cong K^{**}$ defined as $\la b,\,\phi \ra =\la 1,\, \beta \phi \ra$
for all $\phi\in K^*$.
We can compute:
\begin{eqnarray*}
\la b,\,\phi\1 \ra  \phi\2
&=& \la 1,\, \beta \phi\1 \ra \phi\2
     =  \la 1,\, \beta\1\phi\1 \ra \beta\2\phi\2 =\beta\phi, \\
\phi\1 \la b,\,\phi\2 \ra
&=& \phi\1 \la 1,\, \beta \phi\2 \ra
     = \eps\1\phi \la 1,\, \beta\eps\2 \ra =\beta\phi,
\end{eqnarray*}
therefore, $b\in Z(K)$. Also, for all $\phi\in K^*$ we have
\begin{eqnarray*}
\la \eps_s(b),\, \phi \ra
&=& \la b,\, \eps_s(\phi) \ra = \la 1,\, \beta\eps_s(\phi) \ra \\
&=& \la 1,\, \beta\eps\1 \ra \la 1,\, \eps\2\phi \ra =
    \la 1,\, \beta\phi \ra = \la b,\,\phi \ra,
\end{eqnarray*}
therefore $\eps_s(b)=b$ and $b\in K_s$. Thus, $Z(K)\cap K_s \neq
\C\, 1$, so $K$ is not connected.

(ii)$\implies$(i).
If $K$ is not connected, then there exists
$b\in Z(K)\cap K_t,\, b\not\in\C$. Define $\beta\in K^*$
by $\beta : x\mapsto \eps(bx)$.
We have, for all $x\in K$ :
\begin{eqnarray*}
\la \beta,\,\eps_s(x)\ra
&=& \eps(b 1\1)\eps(x 1\2) = \eps(x\eps_t(b)) = \eps(xb), \\
\la \beta,\,\eps_t(x)\ra
&=& \eps(b 1\2)\eps(1\1 x) = \eps(bx) = \eps(xb),
\end{eqnarray*}
from where $\eps_s(\beta)= \beta= \eps_t(\beta)$ and
$K_s^* \cap K_t^* \neq \C\, \eps$.
 
(i)$\implies$(iii).
If there is  a proper subprojection $q$ of $p_\eps$ then from the
formula for $\Delta(p_\eps)$ we get $\eps_s(q)\neq 1$ and
 $\eps_s(q) \in Z(K)$, so $K$ is not connected.
 
(iii)$\implies$(i).
Let $P_\eps$ be the central support of $p_\eps$. It was shown in
\cite{NV} that the quotient map $K \mapsto P_\eps K$
(which is a homomorphism of weak Kac algebras) is one-to-one
on the Cartan subalgebras. Therefore,  $K_s\cap Z(K)$ is contained
in $Z(P_\eps K)$, and $K_s\cap Z(K) =\C$ when $p_\eps$ is minimal.
\end{proof}
\medskip

A very general method of constructing connected weak Kac algebras from 
two-sided crossed-products was introduced in \cite{Nill}. Let $H$ be a
usual finite-dimensional Kac algebra (i.e., finite-dimensional Hopf
$C^*$-algebra) acting on a finite-dimensional $C^*$-algebra $A$ on the
left via $h\otimes a \mapsto (h\lact a)$ and on the right via
$a\otimes h \mapsto (a\ract h)$, $a\in A,\,h\in H$.

\begin{definition}
\label{2-sided}
A two-sided crossed-product $C^*$-algebra $A \rtimes H \ltimes A$
is defined as a vector space $A \otimes H \otimes A$ with multiplication
and involution given by
\begin{eqnarray*}
(b \otimes h \otimes a)(b^\prime \otimes h^\prime \otimes a^\prime) 
&=& b(h\1 \lact b^\prime) \otimes 
    h\2h\1^\prime \otimes (a\ract h\2^\prime) a^\prime  \\
(b \otimes h \otimes a)^* 
&=& (h\1^* \lact b^*) \otimes 
    h\2^* \otimes (a^*\ract h\2^*), 
\end{eqnarray*}
for all $a, a^\prime, b, b^\prime \in A,\, h, h^\prime \in H$.
\end{definition}
\medskip

Let $\{ f_{rs}^\alpha \}$ be a system of matrix units in 
$A = \oplus_\alpha\, M_{m_\alpha}(\C)$ and $\tau$ be the trace of 
the left regular representation of $A$ on itself. Given a
*-anti-automorphism $S_A$ of $A$ define
$$
1\1 \otimes 1\2 = \Sigma_{\alpha rs}\, \frac{1}{m_\alpha}
f_{rs}^\alpha \otimes S_A(f_{sr}^\alpha)
$$
(any projection $P\in A\otimes A$ with the properties
$\mu(\id \otimes S_A^{-1})P = \mu(S_A \otimes \id)P =1$ and
$(\id \otimes \tau)P = (\tau\otimes \id)P =1$, where $\mu$
denotes the multiplication in $A$, has such a form, cf. \cite{NV},
Lemma 2.3.6).

Suppose that the above actions of $H$ on $A$ satisfy
$$
1\1 \otimes (h\lact 1\2) = (1\1\ract h) \otimes 1\2 
\quad \forall h\in H.  
$$

\begin{proposition}[\cite{Nill}]
\label{WKA from 2-sided}
There is a structure of a weak Kac algebra on $K =A \rtimes H \ltimes A$
defined by
\begin{eqnarray*}
\Delta(b \otimes h \otimes a)
&=& (b \otimes h\1 \otimes 1\1) \otimes ((h\2\lact 1\2) \otimes h\3 \otimes a)\\
\eps(b \otimes h \otimes a) &=& \tau( (S_H(h)\lact b)a) \\
S(b \otimes h \otimes a) &=& S_A(a) \otimes S_H(h) \otimes S_A^{-1}(b)
\end{eqnarray*}
\end{proposition}
\medskip

The source and target Cartan subalgebras of $K$ are
$$
K_s =\{1 \otimes 1\otimes a \mid a\in A \}, \qquad
K_t =\{b \otimes 1\otimes 1 \mid b\in A \}.
$$
Clearly, $K_s \cap K_t =\C$, so $K^*$ is connected by 
Proposition~\ref{connectedness}.  It is easy to check that if 
the fixed points algebras
\begin{eqnarray*}
A^H &=& \{ b\in A \mid h\lact b =\eps(h)b\quad \forall h\in H \}, \\
{}^H\!A &=& \{ a\in A \mid a\ract h =\eps(h)a\quad \forall h\in H \}
\end{eqnarray*}
are trivial, then $K_s \cap Z(K) =\C$, and $K$ is biconnected.

In the special case
when $H =\C$ acts trivially on $A$, $K^*$ is isomorphic
to the full matrix algebra $M_d(\C)$, $d=\dim A$.
Such weak Kac algebras were classified in \cite{NV}.
 
\begin{example}
\label{example WKA}
Let $A$ be a Kac subalgebra of $H^*$ and actions of $H$ on $A$
be induced by the dual actions of $H$ on $H^*$ :
$$
h\lact a =a\1 \la h,\, a\2\ra, \qquad
a\ract h = \la h,\, a\1 \ra a\2.
$$
Define $1\1 \otimes 1\2 =\Delta(p)$, where $p$ is the 
Haar projection in $A$. Since $p$ is cocommutative, the condition 
$1\1 \otimes (h\lact 1\2) = (1\1\ract h) \otimes 1\2$
is satisfied for all $h\in H$ and $K= A \rtimes H \ltimes A$ 
is a biconnected weak Kac algebra with $\lambda^{-1} = (\dim H)(\dim A)$.
\end{example}
\medskip

In Section 6 we derive some arithmetic properties of biconnected
weak Kac algebras from the existence of a minimal action of any such
algebra on the hyperfinite II${}_1$ factor.

\end{section}

\begin{section}
{Duality for actions}

In this section $K$ is a weak Kac algebra satisfying the $\lambda$-Markov
condition (e.g., indecomposable) and acting on a $C^*$-algebra $A$.
Left actions are assumed everywhere; the right counterparts of the results
below can be obtained similarly and are left to the reader.

\begin{lemma}
\label{subalgebra}
For all $a\in A$,
$$
(n\lact a) = a(n\lact 1),\,n\in K_s \qquad
(n\lact a) = (n\lact 1)a,\,n\in K_t.
$$
\end{lemma}
\begin{proof}
For all $n\in K_s$ we have
$$
n\lact a = (n\1\lact a)(n\2\lact 1)=(1\1 \lact a)(1\2n\lact 1) = a(n\lact 1),
$$
and similarly the second statement.
\end{proof}

\begin{proposition}
\label{conditional expectation}
The map $E_A : A\rtimes K \to A$, defined as
$$
E_A([a\otimes h]) = a(E_t(h)\lact 1),\qquad a\in A,\, h\in K,
$$
is a faithful conditional expectation. 
If $\{y_\nu\}_{\nu=1\dots n}$ is  a basis for $E_t$ as 
in  Theorem~\ref{Markov criterion}(v), then
$\{[1\otimes y_\nu] \}_{\nu=1\dots n}$ is a basis for $E_A$.
\end{proposition}
\begin{proof}
For all $z \in K_t$ we compute 
\begin{eqnarray*}
E_A([a\otimes zh]) 
&=& a (E_t(zh)\lact 1) = a (zE_t(h)\lact 1) \\
&=& a(z\lact 1)(E_t(h)\lact 1) = E_A([a(z\lact 1)\otimes h]),
\end{eqnarray*}
therefore $E_A$ is well-defined on $A\rtimes K$. Clearly,
$E_A\vert_A =\id_A$. Let us check other properties  
(using Lemma~\ref{subalgebra}) :
\begin{eqnarray*}
E_A([a\otimes 1][b\otimes h][c\otimes 1])
&=& E_A([ab(h\1\lact c)\otimes h\2]) \\
&=& ab(h\1\lact c)(E_t(h\2)\lact 1)
 =  ab(E_t(h)\lact c) \\
&=& ab(E_t(h)\lact 1)c  =  aE_A([b\otimes h])c,
\end{eqnarray*}
for all $a,\,b,\,c\in A$ and $h\in K$, so $E_A$ is a conditional expectation.
We have $h =\Sigma_\nu\, y_\nu E_t(y_\nu^* h) = \Sigma_\nu\, E_t(h y_\nu)
y_\nu^*$ for all $h\in K$  by  Theorem~\ref{Markov criterion}(v), so
\begin{eqnarray*}
[a\otimes h]
&=&  \Sigma_\nu\, [a\otimes E_t(h y_\nu)y_\nu^* ]
=  \Sigma_\nu\, [a (E_t(h y_\nu)\lact 1) \otimes 1] [1\otimes y_\nu^*] \\
&=& \Sigma_\nu\,
[E_A([a\otimes h] [1\otimes y_\nu^*])\otimes 1][1\otimes y_\nu^*],
\end{eqnarray*}
applying the involution we get
$$
[a\otimes h] = \Sigma_\nu\,[1\otimes y_\nu]
[E_A([1\otimes y_\nu^*] [a\otimes h]) \otimes 1], \quad
a\in A,\, h\in K.
$$
Therefore, every $x\in A\rtimes K$ can be written as 
$x = \Sigma_\nu [1\otimes y_\nu] [a_\nu \otimes 1]$ for some $a_\nu,\,
\nu =1\dots n$. Since $E_A([1\otimes y_\nu^*y_\kappa]) =\delta_{\nu\kappa}$,
we have 
$$
E_A(x^*x) = \Sigma_{\nu\kappa}\,E_A([a_\nu^* \otimes 1][1\otimes y_\nu^*])
[1\otimes y_\kappa] [a_\kappa \otimes 1] = \Sigma_\nu\, a_\nu^* a_\nu,
$$
and $x =0$ iff $E_A(x^*x)=0$ iff $a_\nu=0\,(\forall\,\nu)$.
This proves that $E_A$ is  faithful and
$\{[1\otimes y_\nu] \}_{\nu=1\dots n}$ is a basis for $E_A$. 
\end{proof}

\begin{remark}
\label{Index E}
$\Index E_A = \Index E_t =\lambda^{-1}$.
\end{remark}

In what follows we consider $C^*$-algebras $A$, $K$, $K^*$, $A\rtimes K$, 
and $K\rtimes K^*$ as subalgebras of $(A\rtimes K)\rtimes K^*$
in an obvious way with inclusion maps denoted by $i_A$, $i_K$ etc.

\begin{lemma}
\label{exe}
Let $e_A = i_{K^*}(\tau) \in (A\rtimes K)\rtimes K^*$. Then
\begin{enumerate}
\item[(i)] $e_A\, i_{A\smrtimes K}(x)\, e_A = i_A(E_A(x))\,e_A$
for all $x\in A\rtimes K$,
\item[(ii)] the map $A\ni a \mapsto i_A(a)e_A \in (A\rtimes K)\rtimes K^*$ 
is injective.
\end{enumerate}
Moreover, $E_{A\smrtimes K}(e_A) = \lambda$.
\end{lemma}
\begin{proof}
For all $a \in A$, $h\in K$ we compute
\begin{eqnarray*}
e_A\, i_{A\smrtimes K}([a\otimes h])\, e_A
&=& [\tau\1\lact[a\otimes h] \otimes \tau\2\tau] \\
&=& [\tau\1\lact[a\otimes h] \otimes \eps_t(\tau\2)]
    [1_{A\smrtimes K} \otimes \tau] \\
&=& [\tau\lact[a\otimes h] \otimes \eps] [1_{A\smrtimes K} \otimes \tau] 
= i_A(E_A([a\otimes h]))\, e_A,
\end{eqnarray*}
which proves (i). Next, we compute
$$
E_{A\smrtimes K}(i_A(a)e_A) = E_{A\smrtimes K}([a\otimes 1] \otimes \tau]) =
[a\otimes 1] (\lambda\eps \lact[1\otimes 1] ) = \lambda i_A(a),
$$
thus proving that the map $a \mapsto i_A(a)e_A$ is injective. 
Taking $a=1$ in the last formula, we obtain $E_{A\smrtimes K}(e_A) = \lambda$.
\end{proof}

\begin{proposition}
\label{BeB}
$(A\rtimes K)\rtimes K^* = (A\rtimes K)e_A(A\rtimes K).$
\end{proposition}
\begin{proof}
Observe that for all $a \in A$, $g,h\in K$
$$
i_{A\smrtimes K}([a\otimes h])\, e_A\, i_K(g) = i_A(a)\,( i_K(h)\,e_A\,i_K(g) ).
$$
Since $(A\rtimes K)\rtimes K^* = \mbox{span}\{i_A(a)i_{K\smrtimes K^*}(x)\mid
a\in A,\,x\in K\rtimes K^*\}$, it suffices to show that 
$K\rtimes K^* = K e_K K$ (here $e_K= [1_K\otimes\tau] \in K\rtimes K^*)$.

For this purpose, we need to show that every element of $K \rtimes~K^*$
can be written as a linear combination of elements
$ i_K(h)\,e_K\,i_K(g) $, $h,g\in K$.

Let $\{ \phi_{ij}^\gamma\}$ be a system of matrix units in $K^*$.
Since $\tau$ is the normalized Haar projection in $K^*$, we have
$$
\Delta(\tau) = \sum_\gamma\frac{1}{c_\gamma}\,
\sum_{ij} \phi_{ij}^\gamma\otimes S(\phi_{ji}^\gamma),
$$
for some integers $c_\gamma$.
Let $\{ v_{ij}^\gamma\}$ be the system of comatrix units in $K$
dual to $\{ \phi_{ij}^\gamma\}$ : 
$\Delta(v_{ij}^\gamma)=\sum_k v_{ik}^\gamma \otimes v_{kj}^\gamma,\,
\eps(v_{ij}^\gamma)=\delta_{ij}$.

Fix $x\in K$ and let  $h_k = xS(v_{pk}^\gamma),\,g_k = c_\gamma 
v_{kl}^\gamma$ $(x\in K)$ for some $\gamma,\,p,\,l$ $(k=1\dots m_\gamma)$.
Then
\begin{eqnarray*}
\sum_k\,i_K(h_k)\, e_K\, i_K(g_k) 
&=& \sum_{kijm}\, [xS(v_{pk}^\gamma)v_{km}^\gamma \otimes 
     \la \phi_{ij}^\gamma,\, v_{ml}^\gamma \ra S(\phi_{ji}^\gamma) ] \\
&=& \sum_{km}\,[xS(v_{pk}^\gamma)v_{km}^\gamma \otimes S(\phi_{lm}^\gamma)] \\
&=& \sum_{m}\, [x\eps_s(v_{pm}^\gamma) \otimes S(\phi_{lm}^\gamma)] \\
&=& [x \otimes \sum_{m}\,\la \eps\1,\,v_{pm}^\gamma\ra 
    \eps\2 S(\phi_{lm}^\gamma)].
\end{eqnarray*}
Since $x\in K$ is arbitrary, it remains to show that elements of the form
$\psi_{lp}^\gamma =\sum_m\,\la \eps\1,\,v_{pm}^\gamma\ra 
\eps\2 S(\phi_{lm}^\gamma)$ form a linear basis for $K^*$. 
We have
\begin{eqnarray*}
\la S(\psi_{lp}^\gamma),\,v_{pq}^\beta\ra
&=& \sum_m\,\la \phi_{lm}^\gamma \eps\1,\, v_{pq}^\beta \ra
    \la \eps\2,\, S(v_{pm}^\gamma) \ra  \\
&=& \sum_{mj}\, \la \phi_{lm}^\gamma ,\, v_{pj}^\beta \ra
    \la \eps,\, v_{jq}^\beta  S(v_{pm}^\gamma) \ra  
 = \delta_{\gamma\beta} \delta_{lp} 
    \sum_m\,\la \eps,\, v_{mq}^\gamma  S(v_{pm}^\gamma) \ra  \\
&=& \delta_{\gamma\beta} \delta_{lp} \eps(S(v_{pq}^\gamma)) =
\delta_{\gamma\beta} \delta_{lp} \delta_{pq},      
\end{eqnarray*}
therefore, $\psi_{lp}^\gamma = S(\phi_{lp}^\gamma)$.
\end{proof}

\begin{corollary}
\label{isomorphism with basic construction}
$(A\rtimes K)\rtimes K^* \cong \la A\rtimes K,\,e_A\ra $,
i.e., $(A\rtimes K)\rtimes K^*$ is the basic construction for
the conditional expectation $E_A$.
\end{corollary}
\begin{proof}
Propositions~\ref{conditional expectation} and \ref{BeB} show that
$(A\rtimes K)\rtimes K^*$ is generated by $A\rtimes K$ and projection
$e_A$, in the way characterizing the basic construction (see 2.3).
\end{proof}
\medskip

The following result is an analogue of the Takesaki duality theorem
for actions of Kac algebras \cite{ES-Takesaki} and Hopf algebras
\cite{BM}.

\begin{theorem} [Duality for actions]
\label{Takesaki duality}
Let $K$ be a weak Kac algebra satisfying the $\lambda$-Markov condition,
acting on a $C^*$-algebra $A$. Then
$$
(A\rtimes K)\rtimes K^* \cong A \otimes M_n(\C), \qquad where~n=\lambda^{-1}.
$$
\end{theorem}
\begin{proof}
By Proposition~\ref{conditional expectation} there is a basis
for $E_A$, therefore $\la A\rtimes K,\, e_A\ra \cong  A\otimes M_n(\C)$,
and the result follows from Corollary~\ref{isomorphism with basic construction}.
\end{proof}

\begin{lemma}
\label{inclusion into the centralizer}
Let $K$ be a weak Kac algebra acting on the right on a $*$-algebra $A$. 
Then $K_t \subset A^\prime \cap K \ltimes A$. 
\end{lemma}
\begin{proof}
If $z\in K_t$, then
\begin{eqnarray*}
i_A(a)i_K(z)
&=& [ z\1 \otimes (a\ract z\2)] = [z1\1 \otimes (a\ract 1\2)] \\
&=& [ z\otimes a] = i_K(z) i_A(a),
\end{eqnarray*}
thus,  $K_t \subset A^\prime \cap K \ltimes A$.
\end{proof}  
 
\begin{definition}
\label{minimality}
A right action of $K$ on $A$ is { \em minimal} if
$K_t = A^\prime \cap K \ltimes A$.
\end{definition}
 
\end{section}


\begin{section}
{Construction of a minimal action of a biconnected weak Kac algebra
on the hyperfinite II${}_1$ factor.}

In this section we assume that $K$ is a biconnected weak Kac algebra,
in particular it satisfies the $\lambda$-Markov condition for some
$\lambda= n^{-1}$.

\begin{lemma}
\label{iterated crossed products}
Let $K$ act on a finite-dimensional $C^*$-algebra $A$. Suppose that
$\tr$ is a trace on $A\rtimes K$, 
and  $E_A$ from Proposition~\ref{conditional expectation}
is the $\tr$-preserving conditional expectation. 
Then $\tr_1 = tr\circ E_{A\smrtimes K}$ is a trace on
$\la A\rtimes K,\, e_A \ra$, extending $\tr$ and satisfying
$E_{A\smrtimes K}(e_A) = \lambda$. 
In other words, if $\tr$ is a trace on $A\rtimes K$ such that $E_A$
preserves it, then $\tr$ is a $\lambda$-Markov trace for the inclusion
$A \subset A\rtimes K$, and $\tr_1$ is its $\lambda$-Markov extension
to $\la A\rtimes K,\,e_A \ra$.
\end{lemma}    
\begin{proof}
Clearly,  $\tr_1$ is a positive functional on $\la A\rtimes K,\, e_A \ra$
extending $\tr$. Let us show that $\tr_1$ is a trace. By
Lemma~\ref{exe}, $E_{A\smrtimes K}(e_A) = \lambda$, therefore
\begin{eqnarray*}
\tr_1((x_1e_Ay_1)(x_2e_Ay_2)) &=& \tr_1((x_1E_A(y_1x_2)e_Ay_2) =
\lambda\tr(E_A(x_1y_2)E_A(y_1x_2))  \\
&=& \tr_1((x_2e_Ay_2)(x_1e_Ay_1)),
\end{eqnarray*}
for all $x_1,y_1,x_2,y_2 \in A\rtimes K$. 
Since $\la A\rtimes K,\, e_A \ra$ is spanned
by elements of the form $x e_A y,\,(x,y\in A\rtimes K)$ the result
follows from (\cite{JS}, 3.2.5).
\end{proof}

\begin{remark}
In conditions of Lemma~\ref{iterated crossed products}, $e_A$ is the
Jones projection for the inclusion $A\subset A\rtimes K$ with respect to
the Markov trace $\tr$ and 
$E_{A\smrtimes K} : \la A\rtimes K,\, e_A \ra \to  A\rtimes K$ 
is the $\tr$-preserving conditional expectation.
\end{remark}
\medskip

Note that the map $\phi \mapsto (\phi\lact 1)$ gives an isomorphism
between $K_t^*$ and $K_s$ in the crossed product algebra $K\rtimes K^*$.

\begin{proposition}
\label{symmetric commuting square}
Suppose that $K$ is connected and let $\tr$ be the unique Markov trace for 
the inclusion $i_K(K) \subset K\rtimes K^*$. Then
$$
\begin{array}{ccc}
i_K(K) & \subset & K\rtimes K^* \\
\cup   &         & \cup         \\
i_K(K_s)\equiv i_{K^*}({K_t}^*) & \subset & i_{K^*}(K^*) \\
\end{array}
$$
is a symmetric commuting square with respect to $tr$. 
\end{proposition}
\begin{proof}
By Corollary~\ref{tau is a Markov trace}, $\tau$ is a Markov trace
for the inclusion $K_t \subset K$, and $E_t$ is the $\tau$-preserving
conditional expectation.
Since $K\rtimes K^* = (K_t\rtimes K)\rtimes K^*$, it follows from
Lemma~\ref{iterated crossed products} that 
$\tr$ extends $\tau$ and $E_K : K\rtimes K^* \to K$
is the $\tr$-preserving conditional expectation. We have
$$
E_K(i_{K^*}(\phi)) = E_K([1\otimes\phi]) =
i_K(\phi\lact 1) \in i_K(K_s),
$$
for all $\phi\in K^*$.
This proves that the square is commuting. It is symmetric since
$K\rtimes K^* = i_K(K)\,i_{K^*}(K^*)$. 
\end{proof}
\medskip

Corollary~\ref{isomorphism with basic construction}
implies that the sequence
$$
K_t \subset K \subset K\rtimes K^*  \subset
K\rtimes K^* \rtimes K \subset\cdots \subset M
$$
is the  Jones tower for the inclusion $K_t\subset K$ 
on $K$). When $K$ is connected, all the inclusions in this
sequence are connected and the union of these $C^*$-algebras 
admits a unique tracial state, consequently, its von Neumann algebra
completion $M$ with respect to this trace is a copy of the
hyperfinite II${}_1$ factor. Using the standard procedure of iterating
the basic construction we can construct a von Neumann subalgebra $N\subset M$
from the above symmetric commuting square.
 
\begin{proposition}
\label{the lattice}
The lattice of $C^*$-algebras obtained by iterating the basic
construction (in the horisontal direction) for the symmetric
commuting square from Proposition~\ref{symmetric commuting square}
is given by two sequences of alternating crossed product with
$K$ and $K^*$ :
$$
\begin{array}{ccccccccl}
K     & \subset & K\rtimes K^* & \subset & K \rtimes K^* \rtimes K 
      & \subset &    \cdots    & \subset & M                       \\
\cup  &         &     \cup     &         & \cup     
      &         &              &         & \cup                     \\
K_s\equiv K_t^*& \subset &K^* & \subset & K^*\rtimes K 
      & \subset &   \cdots     & \subset &  N, \\
\end{array}
$$
where we identify all $C^*$-subalgebras with their images in $M$.
\end{proposition}
\begin{proof}
Identities $K^*\rtimes K = \la K,\,e_K\ra$,
$K^*\rtimes K \rtimes K^* = \la  K^*\rtimes K,\,e_{K\smrtimes K^*} \ra$
etc. follow immediately from Proposition~\ref{BeB}.
\end{proof}

\begin{proposition}
\label{isomorphism}
There is a *-isomorphism between finite dimensional
$C^*$-algebras 
$$ 
A^r = \underbrace{ K \rtimes K^* \rtimes \cdots \rtimes K \rtimes K^*}_{2r~
factors} \quad  and \quad 
B^r = \underbrace{K \ltimes K^* \ltimes \cdots \ltimes K \ltimes K^*}_{2r~
factors} 
$$
given by the `identity' map
$$
[ h^1 \otimes \phi^1 \otimes \cdots \otimes  h^r \otimes \phi^r]
\mapsto [ h^1 \otimes \phi^1 \otimes \cdots \otimes  h^r \otimes \phi^r],
$$
where $h^i \in K,\, \phi^i \in K^*$.
\end{proposition}
\begin{proof}
By the definition of crossed product, the above algebras are isomorphic to
$$
K \underset{K_t = K_s^*}{\otimes} K^* \underset{K_t^* = K_s}{\otimes}
\cdots \underset{K_t = K_s^*}{\otimes} K^*
$$
as vector spaces. By Theorem~\ref{Takesaki duality} we know that these
algebras are isomorphic to $M_{n^r}(\C)\otimes K_s$, where $n=\lambda^{-1}$.
To see that the `identity' map defines a *-algebra isomorphism, it suffices
to note that
\begin{eqnarray*}
&& [ h^1 \otimes \phi^1 \otimes \cdots \otimes  h^r \otimes \phi^r]
\cdot_{A^r} [g^1 \otimes \psi^1 \otimes \cdots \otimes  g^r \otimes \psi^r] \\
&=& [ h^1 (\phi\1^1\lact g^1) \otimes \phi\2^1 (h\1^2 \lact \psi^1)
 \otimes \cdots \otimes  h\2^r (\phi\1^r\lact g^r) \otimes \phi\2^r\psi^r] \\
&=& [ h^1 g\1^1 \otimes \la \phi\1^1, g\2^1 \ra \phi\2^1 \psi\1^1 \otimes
 \cdots \otimes \la \psi\2^{r-1}, h\1^r \ra h\2^r g\1^r \otimes
 \la \phi\1^r, g\2^r \ra \phi\2^r \psi^r ] \\
&=& [h^1 g\1^1 \otimes (\phi^1 \ract g\2^1)\psi\1^1 \otimes  \cdots \otimes
     (h^r \ract \psi\2^{r-1})g\1^r \otimes (\phi^r \ract g\2^{r})\psi^r]\\
&=& [ h^1 \otimes \phi^1 \otimes \cdots \otimes  h^r \otimes \phi^r]
\cdot_{B^r} [g^1 \otimes \psi^1 \otimes \cdots \otimes  g^r \otimes \psi^r], 
\end{eqnarray*}
for all $h^i,g^i\in K,\, \phi^i,\psi^i\in K^*,\, i=1\dots r$, i.e.
multiplications in $A^r$ and $B^r$ are the same.
\end{proof}

\begin{corollary}
\label{isomorphic lattice}
The lattice of algebras from Proposition~\ref{the lattice} is isomorphic to 
$$
\begin{array}{ccccccccl}
K     & \subset & K\ltimes K^* & \subset & K \ltimes K^* \ltimes K
      & \subset &    \cdots    & \subset & M                       \\
\cup  &         &     \cup     &         & \cup
      &         &              &         & \cup                     \\
K_s   & \subset &K^* & \subset & K^*\ltimes K
      & \subset &   \cdots     & \subset &  N. \\
\end{array}
$$
\end{corollary}
\begin{proof}
Clearly, the isomorphisms constructed in Proposition~\ref{isomorphism}
are compatible with all inclusions of the lattice from
Proposition~\ref{the lattice}.
\end{proof}
\medskip

Our next goal is to show that there is a right
action of $K$ on $N$ such that $ M \cong K \ltimes N$. 

\begin{proposition}
\label{an action}
Let $i_K : h \mapsto [ h \otimes \eps \otimes 1 \otimes \cdots ]$
be the inclusion of $K$ in $M$, $E_N : M \to N$ be the trace
preserving conditional expectation. Then the map
$$
x \ract h = \lambda^{-1}\, E_N(i_K(p_\eps)\, x\, i_K(h)),\quad
x\in N,\,h\in K
$$
defines a right action of $K$ on $N$ such that $M = K \ltimes N$
(cf. \cite{S}, Sect.\ 5).
\end{proposition}
\begin{proof}
There is a right action of $K$ on the *-subalgebra given
by the union of the generating sequence of $C^*$-algebras of $N$:
$$
[ \phi \otimes g \otimes \cdots ] \ract h = 
[ (\phi\ract h) \otimes g \otimes \cdots ], \quad h,\, g \in K,\, 
\phi\in K^*.
$$
We have 
\begin{eqnarray*}
[ \phi \otimes g \otimes \cdots ] \ract h 
&=&  \lambda^{-1} [(\eps\ract E_s(p_\eps))(\phi\ract h) \otimes 
                                                g \otimes \cdots ] \\
&=& \lambda^{-1} E_N( [p_\eps\otimes (\phi\ract h) \otimes
                                                g \otimes \cdots ] ) \\
&=& \lambda^{-1} E_N(i_K(p_\eps)\, [\phi \otimes g\otimes \cdots ]\, i_K(h) ),
\end{eqnarray*}
therefore, the map $x \ract h = \lambda^{-1}\, E_N(i_k(p_\eps) x i_k(h))$
extends the above action to a weakly continuous action of $K$ on $N$.
Clearly, $K \ltimes N = i_k(K)N = M$.
\end{proof}

\begin{corollary} $[M:N] = \lambda^{-1}$.
\end{corollary}
\begin{proof}
Follows from Remark~\ref{Index E} and Proposition 5.1.9 in \cite{JS}.
\end{proof} 
\medskip

Let us compute the higher relative commutants of the inclusion
$N\subset M$.

\begin{lemma}
\label{commutants}
Let $K$ act on the left on a $C^*$-algebra $A$ then
$$
i_{K^*}(K^*)'\cap i_{A\smrtimes K}(A\rtimes K) \cap (A\rtimes K)\rtimes K^*
= i_A(A).
$$
\end{lemma}
\begin{proof}
Let $C=i_{K^*}(K^*)'\cap i_{A\smrtimes K}(A\rtimes K)
\cap (A\rtimes K)\rtimes K^*$ and $x\in C$. 
Recall that $e_A = i_{K^*}(\tau)$. Then $xe_A =e_Axe_A =E_A(x)e_A$
and since the map $A\rtimes K\ni x\mapsto i_{A\smrtimes K}(x)e_A$ is
injective (Lemma~\ref{exe}), 
it follows that $x\in i_A(A)$ and $C\subset i_A(A).$
 
Conversely, for all $a\in A,\,\phi\in K^*$ we have
\begin{eqnarray*}
i_{K^*}(\phi)\,i_A(a)
&=& [1_{A\smrtimes K}\otimes \phi][[a\otimes 1]\otimes \eps]
= [(\phi\1\lact[a\otimes 1]) \otimes \phi\2] \\
&=& [[a\otimes 1](\phi\1\lact[1\otimes 1])\otimes \phi\2]
    = [[a\otimes 1]\otimes \phi] = i_A(a)\,i_{K^*}(\phi),
\end{eqnarray*}
therefore $i_A(A) = C$.
\end{proof}

\begin{proposition}
Let $N\subset M=M_0 \subset M_1 \subset M_2 \cdots$ be
the Jones tower constructed from the inclusion $N\subset M$. Then
\begin{eqnarray*}
N^\prime \cap M_n 
&\cong&  \underbrace{\cdots \ltimes K \ltimes K^*}_{n~factors}\ltimes K_t,
\quad n \geq 0 \\
M^\prime \cap M_n &\cong& 
\underbrace{ \cdots  \ltimes K^* \ltimes K}_{(n-1)~factors} \ltimes K_t^*,
\quad n\geq 1.
\end{eqnarray*}
In particular, the action of $K$ is minimal.
\end{proposition}
\begin{proof}
Iterating the basic construction for the commuting square from 
Proposition~\ref{symmetric commuting square} in the vertical direction
and using Proposition~\ref{isomorphism} we get the lattice
$$
\begin{array}{ccc}
\cdots &   & \cdots \\
\cup &    & \cup \\
K^*\rtimes K & \subset  &  K^*\rtimes K \rtimes K^* \\
\cup &    & \cup \\ 
K  & \subset &  K \rtimes K^* \\
\cup &    & \cup \\
K_t^* \equiv K_s & \subset & K^* .\\
\end{array}
$$
The Ocneanu compactness argument (\cite{JS}, 5.7) and
Lemma~\ref{commutants} imply that
$$  
N^\prime \cap M = K_t,\quad N^\prime \cap M_1 = K^*,\quad
N^\prime \cap M_2 = K \rtimes K^* \dots
$$
Similarly, one computes the relative commutants for $M$.
\end{proof}

\begin{corollary}[\cite{NSzW}]
The inclusion $N\subset M$ is of depth 2.
\end{corollary}
\begin{proof}
We have seen in Section 4 that $K \rtimes K^* \cong K_t \otimes M_n(\C)$,
where $n =\lambda^{-1}$. Therefore, $\dim\,Z(N^\prime \cap M) =
\dim\,Z(N^\prime \cap M_2)$, and $N\subset M$ is of depth 2.
\end{proof}

\begin{corollary}
The $\lambda$-lattice of higher relative commutants \cite{Popa}
of the inclusion $N\subset M$ is given by
$$
\begin{array}{ccccccccccc}
\C & \subset & K_t^* \equiv K_s& \subset& 
K     & \subset & K^*\ltimes K & \subset & K \ltimes K^* \ltimes K   
      & \subset &    \cdots                           \\
   &         & \cup &    &
\cup  &         &     \cup     &         & \cup
      &         &                                   \\
   &         & \C   & \subset &
K_t\equiv K_s^*& \subset &K^* & \subset & K\ltimes K^*
      & \subset &   \cdots     \\
\end{array}
$$
\end{corollary}

\begin{remark}
In a similar way one can construct a left minimal action of a biconnected
weak Kac algebra on the hyperfinite II${}_1$ factor.
\end{remark}

\end{section}


\begin{section}
{Examples of subfactors and arithmetic properties of
biconnected weak Kac algebras}

Let $K$ be a biconnected weak Kac algebra. Recall the notation
$$
K\cong \oplus_{i=1}^N\,M_{d_i}(\C), \qquad
K_s \cong K_t \cong \oplus_{\alpha=1}^L\,M_{m_\alpha}(\C),
$$
from Section 2.1. Let us also denote 
$d=\dim K_s$. We have $\dim K = d\lambda^{-1}$.
\medskip

Reducing the inclusion $N \subset M = K \ltimes N$ constructed in 
Section 5 by a minimal projection $q\in N^\prime \cap M = K_t$
we get an irreducible inclusion $qN \subset qMq$ of hyperfinite
$II_1$ factors with index $[qMq : qN] = \tau(q)^2 \lambda^{-1}$,
where $\tau$ is the normalized trace on $M$ 
( $qN \subset qMq$ is of finite depth \cite{Bisch}, and therefore extremal,
see \cite{Popa1} 1.3.6). 
But $\tau(q) = \frac{m_\alpha}{d}$, when $q\in M_{M_\alpha}(\C)$,
therefore
$$
[qMq : qN] = \frac{m_\alpha^2\,\lambda^{-1}}{d^2}.
$$
Note that since $qN \subset qMq$ has a finite depth, its index is 
an algebraic integer. But by Theorem~\ref{Markov criterion}, $\lambda^{-1}$
is an integer, so $[qMq : qN]$ is rational. Therefore, $[qMq : qN]$
is an integer.

\begin{corollary}
$d^2$ divides $m_\alpha^2\lambda^{-1}$ for all $\alpha$.
\end{corollary}

\begin{corollary}
If $\lambda^{-1} =p$ is a prime, then $K\cong \C\Bbb{Z}_p$.
\end{corollary}
\begin{proof}
By the previous corollary we must have $d=1$, so $\dim K =
d\lambda^{-1} =p$ and the result follows from Corollary~\ref{dim p}.
\end{proof}
\medskip

Next, reducing the inclusion $M \subset M_2$ by a minimal projection
$q$ from the relative commutant
$M^\prime \cap M_2 =K$ we get an irreducible inclusion 
$qM \subset qM_2q$. Clearly, this inclusion depends only on the equivalence
class of $q$, so inclusions of the above type are in one-to-one 
correspondence with irreducible representations of $K$. The index is
$$  
[qM_2q : qM] = \tau(q)^2 [M_2:M] =\tau(q)^2 \lambda^{-2} 
           = (\frac{d_i}{d})^2,
$$
whenever $q\in M_{d_i}(\C)$. Again, the index must be an integer, so we
get the following arithmetic property of biconnected weak Kac algebras.

\begin{corollary}
The dimension of a Cartan subalgebra of $K$ divides the degree of any 
irreducible representation of $K$, i.e. $d$ divides $d_i$ for all $i$.
In particular, $d^2$ divides $\dim K$, and $d$ divides 
$\lambda^{-1} =[M:N]$.
\end{corollary}
\medskip

Finally, let us remark that considering the biconnected weak Kac algebra
$K = H \rtimes H^* \ltimes H$    
constructed from a Kac algebra $H$ as in Example~\ref{example WKA},
we can associate an irreducible subfactor with any 
irreducible representation of $H$ (since we have $K_t =H$ in this case).

\end{section}


\end{document}